\documentclass[11pt]{article}

\usepackage{latexsym}
\usepackage{amssymb}

\newtheorem{thm}{Theorem}
\newtheorem{cor}{Corrolary}

\begin{document}
{
\begin{center}
{\Large\bf
The Nevanlinna-Pick matrix interpolation in the Carath\'eodory class with infinite data both in
the nondegenerate and degenerate cases.}
\end{center}
\begin{center}
{\bf S.M. Zagorodnyuk}
\end{center}

\section{Introduction.}

The classical Nevanlinna-Pick interpolation problem for the case of finite data appeared in
papers~\cite{cit_1000_P}, \cite{cit_2000_N}, and for the case of infinite data in~\cite{cit_3000_N}.
Various matrix and operator-valued generalizations were introduced and investigated by different approaches afterwards.
For a detailed exposition of the subject we refer
to books~\cite{cit_4000_BGR}, \cite{cit_5000_D}, \cite{cit_6000_FF}, \cite{cit_7000_RR},
\cite{cit_8000_GS}.

In this paper we shall analyze the following problem (see Notations below).
Let $\{ z_k \}_{k=0}^\rho$, $z_k\in\mathbb{D}$, be prescribed distinct points: $z_j\not= z_l$, $j\not= l$,
$j,l\in\overline{0,\rho}$; $\rho\in \mathbb{N}\cup \{ \infty \}$. Let
$\{ C_k \}_{k=0}^\rho$, $C_k\in \mathbb{C}_{N\times N}$,
be given. The problem is to find a $\mathbb{C}_{N\times N}$-valued analytic function $T(z)$, $z\in \mathbb{D}$,
which belongs to the Carath\'eodory class $\mathcal{C}_N$, subject to conditions:
\begin{equation}
\label{f1_1}
T(z_k) = C_k,\qquad k=\overline{0,\rho}.
\end{equation}
Here $N\in \mathbb{N}$ and $\rho\in \mathbb{N}\cup \{ \infty \}$, are fixed.
This problem is said to be the Nevanlinna-Pick matrix interpolation problem in the Carath\'eodory class
(with finite or infinite data). The problem is said to be determinate if it has a unique solution.

%\noindent
In 1957, Sz\"okefalvi-Nagy and Koranyi presented their famous pure operator approach to the
operator-valued Nevanlinna-Pick interpolation~\cite{cit_9000_SK},\cite{cit_10000_SK}.
They derived conditions of the solvability for various operator-valued Nevanlinna-Pick problems.
In particular, their results apply to the problem~(\ref{f1_1}) with finite data ($\rho<\infty$).
The latter problem ($\rho<\infty$) was investigated by  Chen and Hu both in the nondegenerate and
degenerate cases using a different method~\cite{cit_11000_ChH}.

%\noindent
The aim of our present investigation is to develop the approach of Sz\"okefalvi-Nagy and Koranyi
to obtain an analytic description of solutions for the problem~(\ref{f1_1}) both in
the nondegenerate and degenerate cases.
In order to obtain an analytic description of solutions, we shall
use important results of Chumakin on generalized resolvents of isometric operators~\cite{cit_12000_Ch},
\cite{cit_13000_Ch}.
A similar approach was recently used in~\cite{cit_14000_Z}, \cite{cit_15000_Z}, \cite{cit_16000_Z}
to treat various matrix moment problems.
Also, the necessary and sufficient conditions for the determinacy
of the problem~(\ref{f1_1}) are obtained. They become especially simple in the case $\rho<\infty$.

{\bf Notations. }
As usual, we denote by $\mathbb{R}, \mathbb{C}, \mathbb{N}, \mathbb{Z}, \mathbb{Z}_+$,
the sets of real numbers, complex numbers, positive integers, integers and non-negative integers,
respectively; $\mathbb{D} = \{ z\in \mathbb{C}:\ |z|<1 \}$, $\mathbb{T} = \{ z\in \mathbb{C}:\ |z|=1 \}$.
The set of all complex vectors of size $N$: $a = (a_0,a_1,\ldots,a_{N-1})$, we denote by
$\mathbb{C}_N$, $N\in \mathbb{N}$.
If $a\in \mathbb{C}_N$, then $a^*$ means its complex conjugate vector.
The set of all complex matrices of size $(N\times N)$ we denote by $\mathbb{C}_{N\times N}$.

\noindent
If $\rho\in \mathbb{Z}_+$, the notation $d\in\overline{0,\rho}$ means that
$d\in \mathbb{Z}_+$, $0\leq d\leq \rho$. The notation $d\in\overline{0,\infty}$ means that
$d\in \mathbb{Z}_+$.

\noindent
Let $M(x)$ be a left-continuous non-decreasing matrix function $M(x) = ( m_{k,l}(x) )_{k,l=0}^{N-1}$
on $[0,2\pi]$, $M(0)=0$, and $\tau_M (x) := \sum_{k=0}^{N-1} m_{k,k} (x)$;
$\Psi(x) = ( dm_{k,l}/ d\tau_M )_{k,l=0}^{N-1}$, $N\in \mathbb{N}$.
By $L^2(M)$ we denote a set (of classes of equivalence)
of $\mathbb{C}_N$-valued functions $f$ on $[0,2\pi]$,
$f = (f_0,f_1,\ldots,f_{N-1})$, such that (see, e.g.,~\cite{cit_17000_MM})
$$ \| f \|^2_{L^2(M)} := \int_0^{2\pi}  f(x) \Psi(x) f^*(x) d\tau_M (x) < \infty. $$
The space $L^2(M)$ is a Hilbert space with a scalar product
$$ ( f,g )_{L^2(M)} := \int_0^{2\pi}  f(x) \Psi(x) g^*(x) d\tau_M (x),\qquad f,g\in L^2(M). $$
We denote $\vec e_m = (\delta_{m,0},\delta_{m,1},...,\delta_{m,N-1})$, $0\leq m\leq N-1$, where
$\delta_{m,j}$ is Kronecker's delta function.

\noindent
By the Carath\'eodory class $\mathcal{C}_N$ we mean a set of all analytic $\mathbb{C}_{N\times N}$-valued functions
$T(z)$ in $\mathbb{D}$ such that $T(z)+T^*(z)\geq 0$.

If H is a Hilbert space then $(\cdot,\cdot)_H$ and $\| \cdot \|_H$ mean
the scalar product and the norm in $H$, respectively.
Indices may be omitted in obvious cases.

\noindent For a linear operator $A$ in $H$, we denote by $D(A)$
its  domain, by $R(A)$ its range, by $\mathop{\rm Ker}\nolimits A$
its null subspace (kernel), and $A^*$ means the adjoint operator
if it exists. If $A$ is invertible then $A^{-1}$ means its
inverse. $\overline{A}$ means the closure of the operator, if the
operator is closable. If $A$ is bounded then $\| A \|$ denotes its
norm.
For a set $M\subseteq H$
we denote by $\overline{M}$ the closure of $M$ in the norm of $H$.
For an arbitrary set of elements $\{ x_n \}_{n\in I}$ in
$H$, we denote by $\mathop{\rm Lin}\nolimits\{ x_n \}_{n\in I}$
the set of all linear combinations of elements $x_n$,
and $\mathop{\rm \overline{span}}\nolimits\{ x_n \}_{n\in I}
:= \overline{ \mathop{\rm Lin}\nolimits\{ x_n \}_{n\in I} }$.
Here $I$ is an arbitrary set of indices.
By $E_H$ we denote the identity operator in $H$, i.e. $E_H x = x$,
$x\in H$. If $H_1$ is a subspace of $H$, then $P_{H_1} =
P_{H_1}^{H}$ is an operator of the orthogonal projection on $H_1$
in $H$.

\section{Descriptions of solutions for the Nevanlinna-Pick problem and the determinacy.}

Let $T(z)$ be a solution of the Nevanlinna-Pick problem~(\ref{f1_1}). As usual, an important role
will be played by the following function (kernel):
\begin{equation}
\label{f2_1}
K(u,v) = \frac{1}{ 2(1-u\overline{v}) } \left( T(u) + T^*(v) \right),\quad u,v\in \mathbb{D}.
\end{equation}
Since $T(z)\in \mathcal{C}_N$, it admits the following representation (e.g.~\cite{cit_18000_B}):
\begin{equation}
\label{f2_2}
T(z) = iT_0 + \int\limits_0^{2\pi} \frac{ e^{it} + z }{ e^{it} - z } dF(t),\quad z\in \mathbb{D},
\end{equation}
where $T_0=T_0^*\in \mathbb{C}_{N\times N}$, and $F(t)$ is a non-decreasing left-continuous
$\mathbb{C}_{N\times N}$-valued function on $[0,2\pi]$. Then
$$ K(u,v) = \frac{1}{ 2(1-u\overline{v}) } \int\limits_0^{2\pi} \left(
\frac{ e^{it} + u }{ e^{it} - u } +
\frac{ e^{-it} + \overline{v} }{ e^{-it} - \overline{v} } \right)
dF(t) $$
\begin{equation}
\label{f2_3}
= \int\limits_0^{2\pi} \frac{ 1 }{ (e^{it} - u) } \frac{ 1 }{ \overline{ (e^{it} - v) } } dF(t),\quad
u,v\in \mathbb{D}.
\end{equation}
Consider the following block matrix (the Pick matrix):
\begin{equation}
\label{f2_4}
P_\rho = ( K(z_k,z_l) )_{k,l=0}^\rho.
\end{equation}
Let
\begin{equation}
\label{f2_5}
f(t) = \sum_{k=0}^d \sum_{m=0}^{N-1} \alpha_{k,m} \frac{1}{ e^{it} - z_k } \vec e_m,\quad
\alpha_{k,m}\in \mathbb{C},\ d\in\overline{0,\rho}.
\end{equation}
We may write
$$ 0 \leq \int_0^{2\pi} f(t) dF(t) f^*(t) $$
$$ = \sum_{k,l=0}^d \sum_{m,n=0}^{N-1} \alpha_{k,m}
\overline{ \alpha_{l,n} } \int_0^{2\pi} \frac{1}{ e^{it} - z_k } \vec e_m dF(t)
\frac{1}{ \overline{(e^{it} - z_l)} } \vec e_n^* $$
$$ = \sum_{k,l=0}^d (\alpha_{k,0},\alpha_{k,1},...,\alpha_{k,N-1})
\int_0^{2\pi} \frac{1}{ e^{it} - z_k } \frac{1}{ \overline{(e^{it} - z_l)} } dF(t)
(\alpha_{l,0},\alpha_{l,1},...,\alpha_{l,N-1})^* $$
$$ = \sum_{k,l=0}^d (\alpha_{k,0},\alpha_{k,1},...,\alpha_{k,N-1})
K(z_k,z_l)
(\alpha_{l,0},\alpha_{l,1},...,\alpha_{l,N-1})^* $$
$$ = \Lambda P_d \Lambda^*, $$
where $\Lambda = (\alpha_{0,0},\alpha_{0,1},...,\alpha_{0,N-1},
\alpha_{1,0},\alpha_{1,1},...,\alpha_{1,N-1},...,
\alpha_{d,0},\alpha_{d,1},...,\alpha_{d,N-1})$. Here we have used the rules for operations on
block matrices.
Therefore
\begin{equation}
\label{f2_6}
P_d \geq 0,\qquad d\in\overline{0,\rho}.
\end{equation}
If $\rho < \infty$ the latter relation means that $P_\rho \geq 0$.

%Thus, condition~(\ref{f2_6}) is necessary for the solvability of the Nevanlinna-Pick problem~(\ref{f1_1}).

Conversely, let the Nevanlinna-Pick problem~(\ref{f1_1}) be given and condition~(\ref{f2_6})
be satisfied. Let
$$ P_\rho = (p_{k,l})_{k,l=0}^{ \rho N + N-1 },\qquad p_{k,l}\in \mathbb{C}. $$
Observe that
\begin{equation}
\label{f2_8}
p_{kN+m,lN+n} = \vec e_m K(z_k,z_l) \vec e_n^*,\qquad k,l\in\overline{0,\rho},\quad 0\leq m,n\leq N-1,
\end{equation}
and
\begin{equation}
\label{f2_9}
K(z_k,z_l) = \frac{1}{ 2 (1 - z_k \overline{z_l}) } \left( C_k + C_l^* \right),\qquad k,l\in\overline{0,\rho}.
\end{equation}
By~\cite[Lemma]{cit_10000_SK} there exist a Hilbert space $H$ and a sequence
$\{ x_k \}_{k=0}^{ \rho N + N-1 }$, $x_k\in H$, such that
\begin{equation}
\label{f2_10}
(x_k,x_l)_H = p_{k,l},\qquad k,l\in\overline{0,\rho N + N - 1},
\end{equation}
and $\mathop{\rm \overline{span}}\nolimits\{ x_k \}_{k=0}^{ \rho N + N-1 } = H$.

Assume that $z_0 = 0$. Set
$$ A_0 \sum_{k=1}^d \sum_{m=0}^{N-1} \alpha_{k,m} x_{kN+m} =
\sum_{k=1}^d \sum_{m=0}^{N-1} \frac{ \alpha_{k,m} }{z_k} ( x_{kN+m} - x_m ), $$
\begin{equation}
\label{f2_11}
\alpha_{k,m}\in \mathbb{C},\quad d\in\overline{0,\rho}.
\end{equation}
Let us check that this definition is correct.
Suppose that
$$ \sum_{k=1}^d \sum_{m=0}^{N-1} \alpha_{k,m} x_{kN+m} =
\sum_{k=1}^d \sum_{m=0}^{N-1} \beta_{k,m} x_{kN+m}, $$
where $\alpha_{k,m},\beta_{k,m}\in \mathbb{C},\quad d\in\overline{0,\rho}$.
Set $\gamma_{k,m} = \alpha_{k,m} - \beta_{k,m}$, $1\leq k\leq d$, $1\leq m\leq N-1$.
We may write
$$ I := \left\| \sum_{k=1}^d \sum_{m=0}^{N-1} \frac{ \gamma_{k,m} }{z_k} ( x_{kN+m} - x_m )
\right\|^2_H $$
$$ = \sum_{k,k'=1}^d \sum_{m,m'=0}^{N-1} \frac{ \gamma_{k,m} }{z_k}
\frac{ \overline{ \gamma_{k',m'} } }{ \overline{ z_{k'} } }
( x_{kN+m} - x_m,  x_{k'N+m'} - x_{m'} )_H. $$
By~(\ref{f2_8}),(\ref{f2_9}),(\ref{f2_10}) we may write
$$ ( x_{kN+m} - x_m,  x_{k'N+m'} - x_{m'} )_H =
( x_{kN+m},  x_{k'N+m'} )_H - ( x_{kN+m},  x_{ m'} )_H $$
$$ - ( x_{m},  x_{k'N+m'} )_H +
( x_{m},  x_{m'} )_H = p_{kN+m,k'N+m'} - p_{kN+m,m'} - p_{m,k'N+m'}$$
$$ + p_{m,m'} = \vec e_m ( K(z_k,z_{k'}) - K(z_k,0) - K(0,z_{k'}) + K(0,0) ) \vec e_{m'} $$
$$ = \vec e_m \left( \frac{1}{ 2( 1 - z_k\overline{ z_{k'} } ) } (C_k + C_{k'}^*)  -
\frac{1}{2}(C_k + C_0^*) - \frac{1}{2}(C_0 + C_{k'}^*) \right. $$
$$ \left. + \frac{1}{2}(C_0 + C_0^*) \right) \vec e_{m'} $$
$$ = z_k\overline{ z_{k'} } \vec e_m \left(
\frac{1}{ 2( 1 - z_k\overline{ z_{k'} } ) } (C_k + C_{k'}^*) \right)
\vec e_{m'} = z_k\overline{ z_{k'} } \vec e_m
K(z_k,z_{k'})
\vec e_{m'} $$
\begin{equation}
\label{f2_12}
= z_k\overline{ z_{k'} } p_{kN+m,k'N+m'} = z_k\overline{ z_{k'} } (x_{kN+m},x_{k'N+m'})_H.
\end{equation}
Therefore
$$ I = \sum_{k,k'=1}^d \sum_{m,m'=0}^{N-1} \gamma_{k,m} \overline{ \gamma_{k',m'} }
(x_{kN+m},x_{k'N+m'})_H $$
$$ = \left( \sum_{k=1}^d \sum_{m=0}^{N-1} \gamma_{k,m} x_{kN+m},
\sum_{k'=1}^d \sum_{m'=0}^{N-1} \gamma_{k',m'} x_{k'N+m'} \right)_H = 0. $$
Consequently, the definition of $A_0$ is correct. Let
$$ x = \sum_{k=1}^d \sum_{m=0}^{N-1} a_{k,m} x_{kN+m},\quad
y = \sum_{k'=1}^d \sum_{m'=0}^{N-1} b_{k',m'} x_{k'N+m'}, $$
where $a_{k,m},b_{k,m}\in \mathbb{C}$, $d\in\overline{0,\rho}$.
By~(\ref{f2_12}) we may write
$$ (A_0 x,A_0 y)_H $$
$$ = \sum_{k,k'=1}^d \sum_{m,m'=0}^{N-1} a_{k,m}\overline{b_{k',m'}}
\frac{1}{z_k \overline{ z_{k'} }}
( x_{kN+m} - x_m, x_{k'N+m'} - x_{m'} )_H $$
$$ = \sum_{k,k'=1}^d \sum_{m,m'=0}^{N-1} a_{k,m}\overline{b_{k',m'}}
( x_{kN+m}, x_{k'N+m'} )_H = (x,y)_H. $$
Therefore $A_0$ is an isometric operator in $H$. Set $A = \overline{A_0}$.
By the definition of $A$ we may write
$$ (E_H - z_k A) x_{kN+m} = x_m,\qquad k\in\overline{0,\rho},\quad 0\leq m\leq N-1; $$
\begin{equation}
\label{f2_13}
x_{kN+m} = (E_H - z_k A)^{-1} x_m,\qquad k\in\overline{0,\rho},\quad 0\leq m\leq N-1.
\end{equation}
Let $\widetilde A\supseteq A$ be a unitary operator in a Hilbert space $\widetilde H\supseteq H$.
Recall that the following operator-valued function (\cite{cit_12000_Ch}, \cite{cit_13000_Ch}):
\begin{equation}
\label{f2_14}
\mathbf{R}_z(A) = P^{\widetilde H}_H ( E_{\widetilde H} - z\widetilde A )^{-1},\qquad z\in
\mathbb{C}\backslash \mathbb{T},
\end{equation}
is said to be a \textbf{generalized resolvent of $A$} (corresponding to $\widetilde A$).
%If $\{ \widetilde B_t \}_{t\in [0,2\pi]}$ is the left-continuous orthogonal resolution of
%unity of $\widetilde A$ then the following operator-valued function:
%\begin{equation}
%\label{f2_15}
%\mathbf{B}_t = P^{\widetilde H}_H \widetilde B_t,\qquad t\in [0,2\pi],
%\end{equation}
%is said to be a \textbf{generalized spectral function of $A$} (corresponding to $\widetilde A$).

\noindent
Set
\begin{equation}
\label{f2_16}
T(z) = i \mathop{\rm Im}\nolimits C_0 +
\left(
( [-E_H + 2 \mathbf{R}_z(A)] x_m, x_l)_H
\right)_{m,l=0}^{N-1},\qquad z\in \mathbb{D}.
\end{equation}
Let us check that $T(z)$ is a solution of the Nevanlinna-Pick problem~(\ref{f1_1}).
In fact, the function $T(z)$ has the following representation:
$$ T(z) = i \mathop{\rm Im}\nolimits C_0 +
\left(
( [-E_{\widetilde H} + 2 (E_{\widetilde H} - z \widetilde A)^{-1}] x_m, x_l)_{\widetilde H}
\right)_{m,l=0}^{N-1} $$
$$ = i \mathop{\rm Im}\nolimits C_0 + \int_0^{2\pi} \left( - 1 + \frac{2}{1-ze^{-i\theta}} \right)
d\left( (\widetilde G_\theta x_m,x_l)_{\widetilde H} \right)_{m,l=0}^{N-1} $$
$$ = i \mathop{\rm Im}\nolimits C_0 + \int_0^{2\pi} \frac{e^{i\theta}+z}{e^{i\theta}-z}
d\left( (\widetilde G_\theta x_m,x_l)_{\widetilde H} \right)_{m,l=0}^{N-1}, $$
where $\{ \widetilde G_\theta \}_{\theta\in [0,2\pi]}$ is the left-continuous orthogonal
resolution of unity of the operator $\widetilde A^{-1}$.
Therefore $T(z)\in \mathcal{C}_{N}$. By~(\ref{f2_13}),(\ref{f2_10}),(\ref{f2_8}),(\ref{f2_9}) we may write
$$ T(z_k) = i \mathop{\rm Im}\nolimits C_0 +
\left(
( [-E_{\widetilde H} + 2 (E_{\widetilde H} - z_k \widetilde A)^{-1}] x_m, x_l)_{\widetilde H}
\right)_{m,l=0}^{N-1} $$
$$ = i \mathop{\rm Im}\nolimits C_0 +
\left(
( -x_m + 2 x_{kN+m}, x_l)_{\widetilde H}
\right)_{m,l=0}^{N-1} $$
$$ = i \mathop{\rm Im}\nolimits C_0 +
\left(
 -p_{m,l} + 2 p_{kN+m,l}
\right)_{m,l=0}^{N-1}
$$
$$ = i \mathop{\rm Im}\nolimits C_0 - K(0,0) + 2K(z_k,0) $$
$$ = \frac{1}{2} (C_0 - C_0^*)  - \frac{1}{2}(C_0+C_0^*) + C_k+C_0^*
= C_k,\qquad k\in\overline{0,\rho}. $$
Thus, $T(z)$ is a solution of the Nevanlinna-Pick problem~(\ref{f1_1}).

Let $\widehat T(z)$ be an arbitrary solution of the Nevanlinna-Pick problem~(\ref{f1_1}).
Let us show that it admits a representation of the form~(\ref{f2_16}) with a generalized resolvent of $A$.
Consider the  space $L^2(F)$, where $F=F(t)$ is taken from the representation~(\ref{f2_2}) for
$\widehat T(z)$.
Let
\begin{equation}
\label{f2_17}
f(t) = \sum_{k=0}^d \sum_{j=0}^{N-1} a_{k,j} \frac{1}{e^{it}-z_k} \vec e_j,
\end{equation}
\begin{equation}
\label{f2_18}
g(t) = \sum_{k'=0}^d \sum_{j'=0}^{N-1} b_{k',j'} \frac{1}{e^{it}-z_{k'}} \vec e_{j'},
\end{equation}
where $a_{k,j},b_{k',j'}\in \mathbb{C}$, $d\in\overline{0,\rho}$.
A set in $L^2(F)$ of all (classes of equivalence of) functions of the form~(\ref{f2_17})
we shall denote by $M^2_0(F)$, and $L^2_0(F)=\overline{M^2_0(F)}$. By~(\ref{f2_3}),(\ref{f2_8}),(\ref{f2_10})
we may  write:
$$ (f(t),g(t))_{ L^2(F) } =
\sum_{k,k'=0}^d \sum_{j,j'=0}^{N-1} a_{k,j} \overline{b_{k',j'}}
\left(
\frac{1}{e^{it}-z_k} \vec e_j, \frac{1}{e^{it}-z_{k'}} \vec e_{j'}
\right)_{L^2(F)} $$
$$ = \sum_{k,k'=0}^d \sum_{j,j'=0}^{N-1} a_{k,j} \overline{b_{k',j'}}
\int_0^{2\pi} \frac{1}{e^{it}-z_k} \frac{1}{ \overline{e^{it}-z_{k'}} } \vec e_j dF(t) \vec e_{j'}^* $$
$$ = \sum_{k,k'=0}^d \sum_{j,j'=0}^{N-1} a_{k,j} \overline{b_{k',j'}}
\vec e_j K(z_k,z_{k'}) \vec e_{j'}^* $$
$$ = \sum_{k,k'=0}^d \sum_{j,j'=0}^{N-1} a_{k,j} \overline{b_{k',j'}}
p_{kN+j,k'N+j'} $$
$$ = \sum_{k,k'=0}^d \sum_{j,j'=0}^{N-1} a_{k,j} \overline{b_{k',j'}}
(x_{kN+j}, x_{k'N+j'})_H $$
\begin{equation}
\label{f2_19}
= \left(
\sum_{k=0}^d \sum_{j=0}^{N-1} a_{k,j} x_{kN+j}, \sum_{k'=0}^d \sum_{j'=0}^{N-1} b_{k',j'} x_{k'N+j'}
\right)_H.
\end{equation}
Consider the following operator:
\begin{equation}
\label{f2_20}
W_0 f(t) = \sum_{k=0}^d \sum_{j=0}^{N-1} a_{k,j} x_{kN+j}.
\end{equation}
Let us check that this operator is correctly defined as an operator from
$M^2_0(F)$ to $H$. Let $f(t)$ and $g(t)$ have the form~(\ref{f2_17}),(\ref{f2_18}).
Suppose that they belong to the same class of equivalence in $L^2(F)$:
$\| f(t)-g(t) \|_{L^2(F)} = 0$.
By~(\ref{f2_19}) we may write
$$ 0 =
\left(
\sum_{k=0}^d \sum_{j=0}^{N-1} (a_{k,j}-b_{k,j}) \frac{1}{e^{it}-z_k} \vec e_j,
\sum_{k'=0}^d \sum_{j'=0}^{N-1} (a_{k',j'}-b_{k',j'}) \frac{1}{e^{it}-z_{k'}} \vec e_{j'}
\right)_{L^2(F)} $$
$$ = \left(
\sum_{k=0}^d \sum_{j=0}^{N-1} (a_{k,j}-b_{k,j}) x_{kN+j},
\sum_{k'=0}^d \sum_{j'=0}^{N-1} (a_{k',j'}-b_{k',j'}) x_{k'N+j'}
\right)_H $$
$$ = \left\| \sum_{k=0}^d \sum_{j=0}^{N-1} a_{k,j} x_{kN+j} -
\sum_{k=0}^d \sum_{j=0}^{N-1} b_{k,j} x_{kN+j} \right\|_H^2. $$
Thus, the operator $W_0$ is defined correctly.
Relation~(\ref{f2_19}) implies that $W_0$ is an isometric operator.
Set $W = \overline{W_0}$. The  operator $W$ is a unitary transformation which maps
$L^2_0(F)$ onto $H$.
Set
$$ L^2_1(F) := L^2(F) \ominus L^2_0(F). $$
The operator
$$ U := W \oplus E_{ L^2_1(F) }, $$
is a unitary transformation which maps $L^2(F)$
onto $H_1 := H\oplus L^2_1(F)$.
Consider the following unitary operator in $L^2(F)$:
$$ Q f(t) = e^{-it} f(t),\qquad f(t) \in L^2(F). $$
Then
$$ \widehat A := U Q U^{-1}, $$
is a unitary operator in $H_1$. Observe that
$$ \widehat A x_{kN+j} = U Q \frac{1}{e^{it}-z_k} \vec{e}_j =
U \frac{1}{e^{it}} \frac{1}{e^{it}-z_k} \vec{e}_j $$
$$ = U \frac{1}{z_k} \left( \frac{1}{e^{it}-z_k} - \frac{1}{e^{it}} \right) \vec{e}_j =
\frac{1}{z_k} (x_{kN+j} - x_j), $$
where $k\in\overline{1,\rho}, 0\leq j\leq N-1$.
Therefore $\widehat A\supseteq A$.
Let $\{ \widehat G_\theta \}_{\theta\in [0,2\pi]}$ is the left-continuous orthogonal
resolution of unity of the operator $\widehat A^{-1}$.
We may write:
$$ \int_0^{2\pi} \frac{1}{1-\zeta e^{it}} d( \widehat G_t x_m,x_l)_{H_1} =
\left(
( E_{H_1} - \zeta \widehat A^{-1}  )^{-1} x_m, x_l
\right)_{H_1} $$
$$ =
\left(
U^{-1} ( E_{H_1} - \zeta \widehat A^{-1} )^{-1} U \frac{1}{e^{it}} \vec e_m, \frac{1}{e^{it}} \vec e_l
\right)_{L^2(F)} $$
\begin{equation}
\label{f2_20_1}
=
\left(
( E_{L^2(F)} - \zeta Q^{-1} )^{-1} \vec e_m, \vec e_l \right)_{ L^2(F) }
=
\int_0^{2\pi} \frac{1}{1-\zeta e^{it}} \vec e_m dF(t) \vec e_l^*,\quad \zeta\in \mathbb{D}.
\end{equation}
From~(\ref{f2_2}) it follows that $T_0 = \mathop{\rm Im}\nolimits \widehat T(0) = \mathop{\rm Im}\nolimits C_0$.
By~(\ref{f2_10}),(\ref{f2_2}),(\ref{f2_9}),(\ref{f2_8}) we have
\begin{equation}
\label{f2_20_2}
\int_0^{2\pi} d(\widehat G_t x_m,x_l)_H = (x_m,x_l)_H = p_{m,l};
\end{equation}
$$ \int_0^{2\pi} \vec e_m dF(t) \vec e_l^* = \vec e_m (\widehat T(0) - i T_0) \vec e_l^* $$
$$ = \vec e_m \left( C_0 - i \frac{C_0-C_0^*}{2i} \right) \vec e_l^*
= \frac{1}{2} \vec e_m (C_0 + C_0^*) \vec e_l^* $$
\begin{equation}
\label{f2_20_3}
= \vec e_m K(0,0) \vec e_l^* = p_{m,l},\qquad 0\leq m,l\leq N-1.
\end{equation}
By~(\ref{f2_20_1}),(\ref{f2_20_2}),(\ref{f2_20_3}) we obtain that
$$ \varphi(\zeta) := \frac{1}{2} \int_0^{2\pi} \frac{1+\zeta e^{it}}{1-\zeta e^{it}} d(\widehat G_t x_m,x_l)_H $$
$$ = -\frac{1}{2} \int_0^{2\pi} d(\widehat G_t x_m,x_l)_H
+ \int_0^{2\pi} \frac{1}{1-\zeta e^{it}} d(\widehat G_t x_m,x_l)_H $$
$$ = -\frac{1}{2} \int_0^{2\pi} \vec e_m dF(t) \vec e_l^* +
\int_0^{2\pi} \frac{1}{1-\zeta e^{it}} \vec e_m dF(t) \vec e_l^* $$
\begin{equation}
\label{f2_20_4}
= \frac{1}{2} \int_0^{2\pi} \frac{1+\zeta e^{it}}{1-\zeta e^{it}} \vec e_m dF(t) \vec e_l^*,\qquad
z\in \mathbb{D}.
\end{equation}
By the inversion formula (\cite[p. 50]{cit_19000_AK}) we conclude that
\begin{equation}
\label{f2_21}
F(t) = \left(
( \widehat G_t x_m,x_l)_{H_1}
\right)_{m,l=0}^{N-1} + c,\qquad t\in [0,2\pi],\ c=\mathrm{const}.
\end{equation}
By~(\ref{f2_2}),(\ref{f2_21}) we may write:
$$ \widehat T(z) = iT_0 + \left( \int_0^{2\pi} \frac{e^{it}+z}{e^{it}-z} d(\widehat G_t x_m,x_l)_{H_1}
\right)_{m,l=0}^{N-1} $$
$$ = iT_0 + \left( \left( \int_0^{2\pi} [-1 + 2(1-ze^{-it})^{-1}] d\widehat G_t x_m,x_l \right)_{H_1}
\right)_{m,l=0}^{N-1} $$
$$ = iT_0 + \left( \left( [ -E_{H_1} + 2(E_{H_1}-z\widehat A)^{-1}] x_m,x_l \right)_{H_1}
\right)_{m,l=0}^{N-1} $$
$$ = i \mathop{\rm Im}\nolimits C_0 + \left( \left( [ -E_{H} + 2\mathbf{R}_z(A)] x_m,x_l \right)_{H}
\right)_{m,l=0}^{N-1}, $$
where $\mathbf{R}_z(A)$ is a generalized resolvent of $A$ (which corresponds to $\widehat A$).
Therefore $\widehat T(z)$ has a representation of the form~(\ref{f2_16}).

\begin{thm}
\label{t2_2}
Let the Nevanlinna-Pick problem~(\ref{f1_1}) with $z_0 = 0$ be given and
condition~(\ref{f2_6}) hold. Let an operator $A=\overline{A_0}$ be constructed for the
Nevanlinna-Pick problem as in~(\ref{f2_11}).
All solutions of the Nevanlinna-Pick problem have the following form:
\begin{equation}
\label{f2_22}
T(z) = i \mathop{\rm Im}\nolimits C_0 +
\left(
( [-E_H + 2 \mathbf{R}_z(A)] x_m, x_l)_H
\right)_{m,l=0}^{N-1},\qquad z\in \mathbb{D},
\end{equation}
where $\mathbf{R}_z(A)$ is a generalized resolvent of $A$.
Conversely, an arbitrary generalized resolvent of $A$
generates by formula~(\ref{f2_22})
a solution of the Nevanlinna-Pick problem~(\ref{f1_1}).

\noindent
Moreover, the correspondence between all generalized resolvents of $A$ and all solutions
of the Nevanlinna-Pick problem is bijective.
\end{thm}
{\bf Proof. }
It remains to prove that different generalized resolvents of $A$ produce different
solutions of the Nevanlinna-Pick problem~(\ref{f1_1}).
Suppose that there exist unitary extensions $\widetilde A_j\supseteq A$ in Hilbert spaces
$\widetilde H_j\supseteq H$, $j=1,2$, such that
\begin{equation}
\label{f2_23}
\mathbf{R}_z^1(A) = P^{\widetilde H_1}_H (E_{\widetilde H_1} - z\widetilde A_1)^{-1} \not=
\mathbf{R}_z^2(A) = P^{\widetilde H_2}_H (E_{\widetilde H_2} - z\widetilde A_2)^{-1},
\end{equation}
$$ ( (-E_H + 2\mathbf{R}_z^1(A)) x_m,x_l )_H = ( (-E_H + 2\mathbf{R}_z^2(A)) x_m,x_l )_H, $$
where $0\leq m,l\leq N-1$, $z\in \mathbb{D}$.
Set $L_N = \mathop{\rm Lin}\nolimits\{ x_k \}_{k=0}^{N-1}$,
$L = \mathop{\rm Lin}\nolimits\{ x_k \}_{k=0}^{\rho N + N-1}$. From the latter relation
by the linearity we get
\begin{equation}
\label{f2_24}
( \mathbf{R}_z^1(A) x,y )_H = ( \mathbf{R}_z^2(A) x,y )_H,\qquad x,y\in L_N,\quad z\in \mathbb{D}.
\end{equation}
Set $R_{j,z} = (E_{\widetilde H_j} - z\widetilde A_j)^{-1}$, $j=1,2$. Observe that
$$ R_{j,z} (E_H - zA) x = (E_{\widetilde H_j} - z\widetilde A_j)^{-1}
(E_{\widetilde H_j} - z\widetilde A_j) x = x,\qquad x\in D(A); $$
$$ R_{1,z} h = R_{2,z} h\in H,\qquad h\in (E_H-zA)D(A); $$
\begin{equation}
\label{f2_25}
\mathbf{R}_z^1(A) h = \mathbf{R}_z^2(A) h,\qquad h\in (E_H-zA)D(A),\quad z\in \mathbb{D}.
\end{equation}
Choose arbitrary $z\in \mathbb{D}\backslash\{ 0 \}$, $x\in H$, $h\in (E_H- \frac{1}{\overline{z}} A)D(A)$.
We may write
$$ (\mathbf{R}_z^j(A) x,h)_H = (R_{j,z} x,h)_{H_j} = (x, R_{j,z}^* h)_{H_j}
= (x, (E_{H_j} - R_{ j, \frac{1}{\overline{z}} }) h)_{H_j} $$
$$ = (x,h)_H - (x, \mathbf{R}_{ \frac{1}{\overline{z}} }^j(A) h)_H; $$
\begin{equation}
\label{f2_26}
(\mathbf{R}_z^1(A) x,h)_H = (\mathbf{R}_z^2(A) x,h)_H,\quad x\in H,\ h\in (E_H-\frac{1}{\overline{z}} A)D(A),\
z\in \mathbb{D}\backslash\{ 0 \}.
\end{equation}
Choose arbitrary $x\in L$, $z\in \mathbb{C}\backslash \mathbb{T}$:
$z\not= z_k$, $k=\overline{1,\rho}$. Let us check that there exists the following representation:
\begin{equation}
\label{f2_27}
x = x_0 + x_1,\qquad x_0\in L_N,\quad x_1\in (E_H-zA)D(A).
\end{equation}
Here $x_0$,$x_1$ may depend on the choice of $z$ and $x$.

\noindent
Let $x = \sum_{k=0}^d \sum_{m=0}^{N-1} \alpha_{k,m} x_{kN+m}$, $d\in\overline{0,\rho}$.
Set
$$ \beta_{k,m} = \frac{1}{1-\frac{z}{z_k}} \alpha_{k,m},\qquad 1\leq k\leq d,\quad 0\leq m\leq N-1; $$
$$ h = \sum_{k=1}^d \sum_{m=0}^{N-1} \beta_{k,m} x_{kN+m} \in D(A). $$
Then
$$ x_1 := (E_H-zA) h = h - zAh $$
$$ = \sum_{k=1}^d \sum_{m=0}^{N-1} \beta_{k,m} x_{kN+m}
- z \sum_{k=1}^d \sum_{m=0}^{N-1} \beta_{k,m} \frac{1}{z_k} ( x_{kN+m} - x_m ) $$
$$ = \sum_{k=1}^d \sum_{m=0}^{N-1} \alpha_{k,m} x_{kN+m} + z \sum_{k=1}^d \sum_{m=0}^{N-1} \beta_{k,m}
\frac{1}{z_k} x_m $$
$$ = x - \sum_{m=0}^{N-1} \alpha_{0,m} x_{m} + z \sum_{k=1}^d \sum_{m=0}^{N-1} \beta_{k,m} \frac{1}{z_k} x_m
= x - x_0, $$
where $x_0 := \sum_{m=0}^{N-1} \alpha_{0,m} x_{m} - z \sum_{k=1}^d \sum_{m=0}^{N-1} \beta_{k,m} \frac{1}{z_k}
x_m \in L_N$.
Therefore relation~(\ref{f2_27}) holds.

\noindent
Choose an arbitrary $h\in L$, $z\in \mathbb{D}\backslash\{ 0 \}$. By~(\ref{f2_27}) we may write:
$$ h = h_0 + h_1,\qquad h_0\in L_N,\quad h_1\in (E_H-\frac{1}{\overline{z}}A) D(A). $$
Choose an arbitrary $x\in L_N$. By~(\ref{f2_26}),(\ref{f2_24}) we may write
$$ (\mathbf{R}_z^1(A) x,h)_H = (\mathbf{R}_z^1(A) x,h_0)_H + (\mathbf{R}_z^1(A) x,h_1)_H $$
$$ = (\mathbf{R}_z^2(A) x,h_0)_H + (\mathbf{R}_z^2(A) x,h_1)_H = (\mathbf{R}_z^2(A) x,h)_H. $$
Therefore
\begin{equation}
\label{f2_28}
\mathbf{R}_z^1(A) x = \mathbf{R}_z^2(A) x,\qquad x\in L_N,\quad
z\in \mathbb{D}.
\end{equation}
Choose an arbitrary $g\in L$, $z\in \mathbb{D}:\ x\not= z_k$, $k=\overline{1,\rho}$. By~(\ref{f2_27}) we may write:
$$ g = g_0 + g_1,\qquad g_0\in L_N,\quad g_1\in (E_H-zA) D(A). $$
By~(\ref{f2_25}),(\ref{f2_28}) we get
$$ \mathbf{R}_z^1(A) g = \mathbf{R}_z^1(A) g_0 + \mathbf{R}_z^1(A) g_1
= \mathbf{R}_z^2(A) g_0 + \mathbf{R}_z^2 (A) g_1 = \mathbf{R}_z^2(A) g; $$
$$ \mathbf{R}_z^1(A) = \mathbf{R}_z^2(A),\qquad z\in \mathbb{D}. $$
Since $(\mathbf{R}_z^j(A))^* = E_{H_j} - \mathbf{R}_{ \frac{1}{\overline{z}} }^j(A)$, $j=1,2$, $z\in \mathbb{C}$:
$|z|\not=1$, $z\not=0$ (\cite{cit_13000_Ch}), we conclude that
$\mathbf{R}_z^1(A) = \mathbf{R}_z^2(A)$. We obtained a contradiction with~(\ref{f2_23}).
$\Box$

We shall use the following important result:
\begin{thm}~\cite[Theorem 3]{cit_13000_Ch}
\label{t2_3}
An arbitrary generalized resolvent $\mathbf{R}_\zeta$ of a closed isometric operator $U$
in a Hilbert space $H$ has the following representation:
\begin{equation}
\label{f2_41}
\mathbf R_{\zeta} =
\left[
E - \zeta ( U \oplus \Phi_\zeta )
\right]^{-1},\qquad
\zeta\in \mathbb{D}.
\end{equation}
Here $\Phi_\zeta$ is an analytic in $\mathbb{D}$ operator-valued function which values are
linear contractions (i.e. $\| \Phi_\zeta \| \leq 1$) from $H\ominus D(U)$ into $H\ominus R(U)$.

\noindent
Conversely, each analytic in $\mathbb{D}$ operator-valued function with above properties
generates by relation~(\ref{f2_41}) a generalized resolvent $\mathbf R_{\zeta}$ of $U$.
\end{thm}
Observe that relation~(\ref{f2_41}) also shows that different analytic in $\mathbb{D}$ operator-valued functions
with above properties generate different generalized resolvents of $U$.

Comparing the last two theorems we obtain the following result.
\begin{thm}
\label{t2_4}
Let the Nevanlinna-Pick problem~(\ref{f1_1}) with $z_0=0$ be given and
condition~(\ref{f2_6}) hold. Let an operator $A=\overline{A_0}$ be constructed for the
Nevanlinna-Pick problem as in~(\ref{f2_11}).
All solutions of the Nevanlinna-Pick problem have the following form
\begin{equation}
\label{f2_43}
T(z) = i \mathop{\rm Im}\nolimits C_0 +
\left(
( [-E_H + 2 \left[E - z ( A \oplus \Phi_z ) \right]^{-1}] x_m, x_l)_H
\right)_{m,l=0}^{N-1},\qquad z\in \mathbb{D},
\end{equation}
where $\Phi_z$ is an analytic in $\mathbb{D}$ operator-valued function which values are
linear contractions from $H\ominus D(A)$ into $H\ominus R(A)$.

Conversely, each analytic in $\mathbb{D}$ operator-valued function with above properties
generates by relation~(\ref{f2_43}) a solution of the Nevanlinna-Pick problem~(\ref{f1_1}).

Moreover, the correspondence between all
analytic in $\mathbb{D}$ operator-valued functions with above properties
and all solutions
of the Nevanlinna-Pick problem~(\ref{f1_1}) is bijective.
\end{thm}
{\bf Proof. } The proof is obvious.
$\Box$

\begin{cor}
\label{c2_1}
Let the Nevanlinna-Pick problem~(\ref{f1_1}) with $z_0=0$ be given and
condition~(\ref{f2_6}) hold. Let an operator $A=\overline{A_0}$ be constructed for the
Nevanlinna-Pick problem as in~(\ref{f2_11}).
The Nevanlinna-Pick problem~(\ref{f1_1}) is determinate if and only if
at least one of the defect numbers of $A$ is equal to zero.
\end{cor}
{\bf Proof. } If one of the defect numbers of $A$ is equal to zero, then there exists
the unique function $\Phi_z\equiv 0$ of the required class.
On the other hand, if both defect numbers of $A$ are non-zero, then besides
$\Phi_z\equiv 0$ there exist non-zero suitable constant functions $\Phi_z$.
$\Box$

Conditions for the determinacy of the Nevanlinna-Pick problem~(\ref{f1_1}) become especially simple
in the case $\rho<\infty$.
\begin{thm}
\label{t2_5}
Let the Nevanlinna-Pick problem~(\ref{f1_1}) with $z_0=0$, $\rho<\infty$, be given and
condition~(\ref{f2_6}) hold.
The Nevanlinna-Pick problem~(\ref{f1_1}) is determinate if and only if
at least one of the following conditions hold:
\begin{itemize}
\item[(A)]
For each $k\in\overline{0,N-1}$ the following linear system of equations with unknowns $a_{k,j}$ has a solution:
\begin{equation}
\label{f2_44}
\sum_{j=N}^{\rho N + N-1} a_{k,j} p_{j,l} = p_{k,l},\qquad l=0,1,...,\rho N + N -1.
\end{equation}

\item[(B)]
For each $j\in\overline{0,N-1}$ the following linear system of equations with unknowns $b_{j,kN+m}$ has a solution:
\begin{equation}
\label{f2_45}
\sum_{k=1}^\rho \sum_{m=0}^{N-1} b_{j,kN+m} (p_{kN+m,l}-p_{m,l}) = p_{j,l},\qquad l=0,1,...,\rho N + N -1.
\end{equation}
\end{itemize}
Here $p_{i,j}$ are taken from relation~(\ref{f2_8}).
\end{thm}
{\bf Proof. } Let an operator $A=\overline{A_0}$ be constructed for the
Nevanlinna-Pick problem as in~(\ref{f2_11}).
By~(\ref{f2_10}) we see that condition~(A) is equivalent to the following relation:
$$ x_k = \sum_{j=N}^{\rho N + N-1} a_{k,j} x_j,\qquad k=0,1,...,N-1. $$
The latter relation is equivalent to the condition $D(A)=H$.

\noindent
On the other hand, condition~(B) is equivalent to the following relation:
$$ x_j = \sum_{k=1}^\rho \sum_{m=0}^{N-1} b_{j,kN+m} (x_{kN+m,l}-x_{m}),\qquad j=0,1,...,N-1. $$
The latter relation is equivalent to the condition $R(A)=H$.

\noindent
It remains to apply Corollary~\ref{c2_1} to complete the  proof.
$\Box$

\noindent
{\bf Remark. } As it was noted in~\cite{cit_10000_SK}, condition $z_0=0$ is not restrictive.
If $z_0\not= 0$, one may consider a fractional linear transformation
$$ u = u(z) = \frac{z-z_0}{1-\overline{z_0}z}, $$
and seek for a $\mathbb{C}_{N\times N}$-valued function $R(u)$ in $\mathbb{D}$,
which belongs to $\mathcal{C}_N$, subject to
\begin{equation}
\label{f2_46}
R(u_k) = C_k,\qquad k\in\overline{0,\rho},
\end{equation}
where $u_k := u(z_k)$.
It is easy to see that the following relation:
\begin{equation}
\label{f2_47}
T(z) = R(u(z)),\qquad z\in \mathbb{D},
\end{equation}
establishes a bijective correspondence between all solutions of~(\ref{f2_46}) and all solutions of~(\ref{f1_1}).

\begin{center}
{\large\bf The Nevanlinna-Pick matrix interpolation in the Carath\'eodory class with infinite data both in
the nondegenerate and degenerate cases.}
\end{center}
\begin{center}
{\bf S.M. Zagorodnyuk}
\end{center}

In this paper we study the Nevanlinna-Pick matrix interpolation problem in the Carath\'eodory class
with infinite data (both in the nondegenerate and degenerate cases).
We develop the Sz\"okefalvi-Nagy and Kor\'anyi operator approach to obtain an analytic description
of all solutions of the problem. Simple necessary and sufficient conditions for the determinacy
of the problem are given.

}
\end{document}